\documentclass{amsart}
\usepackage{amsfonts}
\usepackage{amssymb}
\usepackage{amsmath}
\usepackage[dvips]{epsfig,graphics}
\usepackage{latexsym, amsfonts, longtable, tabularx}
\usepackage{epsfig}
\usepackage{psfrag}
\usepackage{amscd}
\usepackage[all]{xy}

\newtheorem{thm}{Theorem}
\newtheorem{prop}{Proposition}
\newtheorem{lemma}{Lemma}
\newtheorem{cor}[thm]{Corollary}
\newtheorem{conj}{Conjecture}
\newtheorem{defn}{Definition}
\newtheorem{rmk}{Remark}

\newtheorem{question}{Question}
\newcommand{\Z}{\mathbb{Z}}

\renewcommand{\c}{\cite}

\newcommand{\pf}{{\em Proof: \quad }}
\newcommand{\done}{\hfill $\blacksquare$}
\newcommand{\kb}[1]{\ensuremath{\langle #1 \rangle}}

\newcommand{\del}{\partial}
\newcommand{\U}{{\widetilde U}}

\newcommand{\UC}{\mathcal{UC}}
\newcommand{\C}{\mathcal{C}}

\begin{document}

\title{On mutation and Khovanov homology}
\author{Abhijit Champanerkar}
\address{Department of Mathematics, College of Staten Island, City University of New York}
\email{abhijit@math.csi.cuny.edu}
\thanks{The first author is supported by NSF grant DMS-0455978}

\author{Ilya Kofman}
\address{Department of Mathematics, College of Staten Island, City University of New York}
\email{ikofman@math.csi.cuny.edu}
\thanks{The second author is supported by NSF grant DMS-0456227 and a PSC-CUNY grant}

\date{June 25, 2008}

\begin{abstract}
\noindent
It is conjectured that the Khovanov homology of a knot is invariant
under mutation. In this paper, we review the spanning tree complex for
Khovanov homology, and reformulate this conjecture using a matroid
obtained from the Tait graph (checkerboard graph) $G$ of a knot
diagram $K$.  The spanning trees of $G$ provide a filtration and a
spectral sequence that converges to the reduced Khovanov homology of
$K$.  We show that the $E_2$--term of this spectral sequence is a
matroid invariant and hence invariant under mutation.
\end{abstract}
\maketitle

\centerline{\em In memory of Xiao-Song Lin}

\section{Introduction}
For any diagram of an oriented link $L$, Khovanov \c{Khovanov}
constructed bigraded abelian groups $H^{i,j}(L),$ whose
bigraded Euler characteristic gives the Jones polynomial $V_L(t)$: 
$$\chi(H^{*,*})= \sum_{i,j}(-1)^i q^j {\rm rank}(H^{i,j})= (q+q^{-1})V_L(q^2) $$ 
For knots (or links with a marked component), Khovanov also defined reduced homology groups
$\widetilde{H}^{i,j}(L)$ whose bigraded Euler characteristic is
$q^{-1}V_L(q^2)$ \c{KhPatterns}.  

Since the introduction of Khovanov homology in \c{Khovanov}, the
theory has been developed and generalized far beyond the Jones polynomial (see e.g. \c{Kh_ICM2006} and references therein), 
and beyond classical links to objects like graphs and ribbon graphs (see e.g. \c{dkh, Rong_HG, LoeblMoffatt}).

However, just as the original Jones polynomial eludes a topological
interpretation in terms of the knot complement, classical Khovanov
homology remains mysterious. The following questions are open for 
this invariant:
\begin{itemize}
\item Does any non-trivial knot have trivial Khovanov homology?
\item Which knots have ``thin'' Khovanov homology (supported on two diagonals)?
\item What is the Khovanov homology of $(p,q)$--torus knots?
\item Is Khovanov homology invariant under Conway mutation of knots?
\end{itemize}
Our purpose here is to present ideas and results that we hope will be useful to tackle the last question.
It is conjectured that the Khovanov homology of a knot is invariant
under mutation (see \c{BNweb, Wehrli2}, and see \c{Wehrli3} for a
recent proof over $\Z/2\Z$).

\subsection*{Spanning trees and Khovanov homology} 

There is a $1$-$1$ correspondence between connected link diagrams $D$
and connected planar graphs $G$ with signed edges.  $G$, called the
{\em Tait graph} of $D$, is obtained by checkerboard coloring complementary
regions of $D$, assigning a vertex to every shaded region, an edge to
every crossing and a $\pm$ sign to every edge as follows:
$$ \includegraphics[height=1cm]{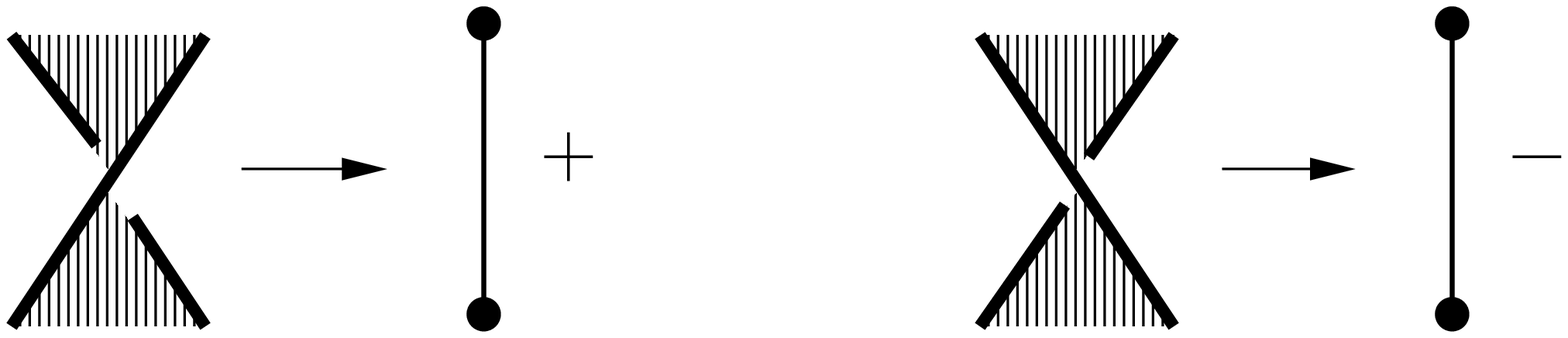}$$
The signs are all equal if and only if $D$ is alternating.  

In Section \ref{Sec_STC} below, we express Khovanov homology using
generators that correspond to spanning trees of $G$.  More accurately,
with a fixed edge order on $G$, the construction relies on {\em
  activity words} $W(T)$ for each spanning tree $T$.  Before diving
into notation, it seems worthwhile to motivate this approach.

We give three motivating reasons to consider Khovanov
homology using the spanning trees of the Tait graph.
First, Thistlethwaite \c{thistlethwaite} gave an expansion of the Jones
polynomial $V_L(t)$ in terms of spanning trees of any Tait graph
$G(L)$.  Every spanning tree contributes a monomial to the Jones
polynomial.  For non-alternating knots, these monomials may cancel
with each other, but for alternating knots, such cancelations do not
occur.  Thus, for alternating knots, the number of spanning trees is
exactly the $L^1$-norm of Jones coefficients, and the span of $V_L(t)$
is maximal, equal to the crossing number.
The bigraded spanning tree complex described below provides an
explicit distribution of spanning trees, which is at most
$(k+1)$--thick for links that become alternating after $k$ crossing
changes (see \c{KHshort}).  It also provides a tool to study
particular Jones coefficients.  For example, if we change a crossing
in an alternating knot diagram $D$, the span and $L^1$-norm of Jones
coefficients strictly decrease.  In the spanning tree complex, we can
see how the gradings change to make certain spanning trees cancel in
the Euler characteristic.  

A second reason is given by the important and closely related example
of knot Floer homology.  The two knot homology theories are compared
in detail in \c{Rasmussen}.  The complex for knot Floer homology in
\c{OZ1} has generators that correspond to spanning trees, but no
combinatorial differential is known.  The more recent complex in
\c{MOS} is completely combinatorial but has far more generators, so it
is quasi-isomorphic (and possibly retracts) to a combinatorial complex
generated by spanning trees.  The situation for Khovanov homology is
similar: The Khovanov complex retracts to a complex generated by
spanning trees of $G$ (Theorem \ref{KHthm}), but it remains an open
question whether the differential on the spanning tree complex can be
expressed entirely in terms of the combinatorics (activity words) of
spanning trees.

The third reason is discussed in Section \ref{Sec_MM}, where we show
that the conjectured dependence of the differential on activity words
is closely related to the mutation invariance of Khovanov homology.
This appears to be a promising approach to prove that Khovanov
homology is invariant under component-preserving mutation of links.

In Section \ref{Sec_STC}, we review the construction of the spanning
tree complex $\C(K)$ given in \c{KHshort}, the spanning tree
filtration and the associated spectral sequence that converges to
$\widetilde{H}(K)$. 
In Sections \ref{Sec_WD} and \ref{Sec_SS}, we prove new results that
show direct incidences and the $E_2$--term of this spectral sequence
are determined by activity words. Material in Section \ref{Sec_higher}
also has not been previously published.

In Section \ref{Sec_MM}, we show that the mutation invariance of any
knot invariant can be expressed in terms of the colored cycle matroid
$M(K)$, obtained from the Tait graph $G$ of a knot diagram $K$.  In
particular, the reduced Khovanov homology $\widetilde{H}(K)$ is
invariant under mutation if and only if the spanning tree complex
$\C(K)$ is determined by $M(K)$ up to quasi-isomorphism.  As a partial
step, the $E_2$--term mentioned above is determined by $M(K)$ and
hence invariant under mutation.  In Section \ref{KHmatroids}, we
discuss an approach to prove mutation-invariance of Khovanov homology.

\section{Spanning tree complex} \label{Sec_STC} 

In \c{KHshort}, for any connected link diagram $D$, we defined the
spanning tree complex $\C(D)=\{\C_v^u(D), \del\}$, whose generators
correspond to spanning trees $T$ of $G$.  In this section, we review
the main ideas and related notation, which will be used later.

\subsection{Activity words and twisted unknots} \label{Sec_awtu}
Fix an order on the edges of $G$.  For every spanning tree $T$ of $G$,
each edge $e\in G$ has an activity with respect to $T$, as follows.
If $e \in T$, $\mathit{cut(T,e)}$ is the set of edges that connect
$T\setminus e$.  If $f \notin T$, $\mathit{cyc(T,f)}$ is the set of
edges in the unique cycle of $T \cup f$.  Note $f \in cut(T,e)$ if and
only if $e \in cyc(T,f)$.  An edge $e \in T$ (resp. $e \notin T$) is
{\em live} if it is the lowest edge in its cut (resp. cycle), and
otherwise it is {\em dead}.

For any spanning tree $T$ of $G$, the {\em activity word} $W(T)$ gives
the activity of each edge of $G$ with respect to $T$.  The letters of
$W(T)$ are as follows: $L,\ D,\ \ell,\ d$ denote a positive edge that
is live in $T$, dead in $T$, live in $G-T$, dead in $G-T$,
respectively; $\bar{L},\ \bar{D},\ \bar{\ell},\ \bar{d}$ denote
activities for a negative edge.  Note that $T$ is given by the capital letters of $W(T)$.

Thistlethwaite assigned a monomial $\mu(T)$ to each $T$ as follows:
$$ L^p D^q \ell^r d^s \bar{L}^x \bar{D}^y \bar{\ell}^z \bar{d}^w \quad \Rightarrow \quad \mu(T)= (-1)^{p+r+x+z}A^{-3p+q+3r-s+3x-y-3z+w} $$
\begin{thm}\c{thistlethwaite}\label{ThistThm}
  Let $G$ be the Tait graph of any connected link diagram $D$ with any
  order on its edges.  Let $\kb{D}$ denote the Kauffman bracket
  polynomial of $D$. Summing over all spanning trees $T$ of $G$,
  $\kb{D}= \sum_T \mu(T)$.
\end{thm}

The activity word $W(T)$ contains much more information than just
$\mu(T)$.  A \emph{twisted unknot} $U$ is a diagram of the unknot
obtained from the round unknot using only Reidemeister I moves.
$W(T)$ determines a twisted unknot $U(T)$ by changing the crossings of
$D$ according to Table \ref{Table1} for dead edges, and leaving the
crossings unchanged for live edges (Lemma 1 \c{KHshort}).  In Table
\ref{Table1}, the sign of the crossing in $U(T)$ is indicated for
unsmoothed crossings, and Kauffman state markers are indicated for
smoothed crossings.
\begin{table}[h]
\caption{Activity word for a spanning tree determines a twisted unknot}
\label{Table1}
\begin{center}
\begin{tabular}{cc|cc|cc|cc}
$L$ & $D$ & $\ell$ & $d$ & $\bar{L}$ & $\bar{D}$ &
$\bar{\ell}$ & $\bar{d}$ \\
\hline
$-$ & $A$ & $+$ & $B$ & $+$ & $B$ & $-$ & $A$ \\
\includegraphics[height=0.5cm]{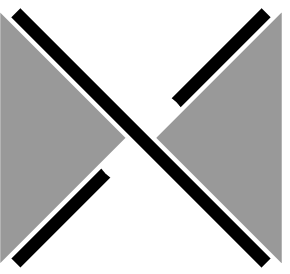} & \includegraphics[height=0.5cm]{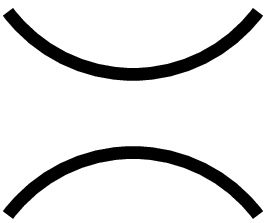} &
\includegraphics[height=0.5cm]{poscross.eps} & \includegraphics[height=0.5cm]{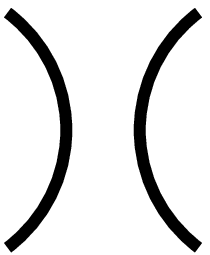} &
\includegraphics[height=0.5cm]{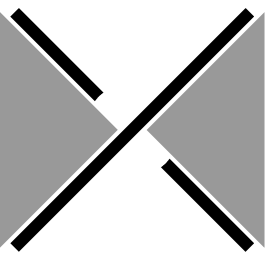} & \includegraphics[height=0.5cm]{posA.eps} &
\includegraphics[height=0.5cm]{negcross.eps} & \includegraphics[height=0.5cm]{posB.eps} \\
\end{tabular}
\end{center}
\end{table}

\begin{figure}
  \begin{center}
    \psfrag{pedge}{positive edge}
    \psfrag{nedge}{negative edge}
    \psfrag{fig8}{$K=$ Figure-8 knot}
    \psfrag{chcolor}{Colored diagram $D$}
    \psfrag{graph}{Signed graph for $D$}
    \psfrag{e1}{\footnotesize{$1$}}
    \psfrag{e2}{\footnotesize{$2$}} \psfrag{e3}{\footnotesize{$3$}}
    \psfrag{e4}{\footnotesize{$4$}}
    \psfrag{sp}{\footnotesize{Spanning}} \psfrag{tr}{\footnotesize{trees}}
    \psfrag{act}{\footnotesize{Activities}}
    \psfrag{wt}{\footnotesize{Weights}}
    \psfrag{a1}{\footnotesize{$LLdd$}} \psfrag{a2}{\footnotesize{$LdDd$}}
    \psfrag{a3}{\footnotesize{$\ell DDd$}}\psfrag{a4}{\footnotesize{$\ell LdD$}}
    \psfrag{a5}{\footnotesize{$\ell \ell DD$}}
    \psfrag{w1}{\footnotesize{$A^{-8}$}} \psfrag{w2}{\footnotesize{$-A^{-4}$}}
    \psfrag{w3}{\footnotesize{$-A^{4}$}} \psfrag{w4}{\footnotesize{$1$}}
    \psfrag{w5}{\footnotesize{$A^{8}$}}
    \psfrag{T1}{\footnotesize{$T_1$}} \psfrag{T2}{\footnotesize{$T_2$}}
    \psfrag{T3}{\footnotesize{$T_3$}} \psfrag{T4}{\footnotesize{$T_4$}}
    \psfrag{T5}{\footnotesize{$T_5$}}
    \psfrag{gamma}{$\langle D \rangle = A^{-8}-A^{-4}+1-A^{4}+A^{8}$}
    \psfrag{writhe}{writhe $w(D)=0$}
    \psfrag{jones}{$V_D(t)=t^{-2}-t^{-1}+1-t+t^2$}
    \includegraphics[width=5in]{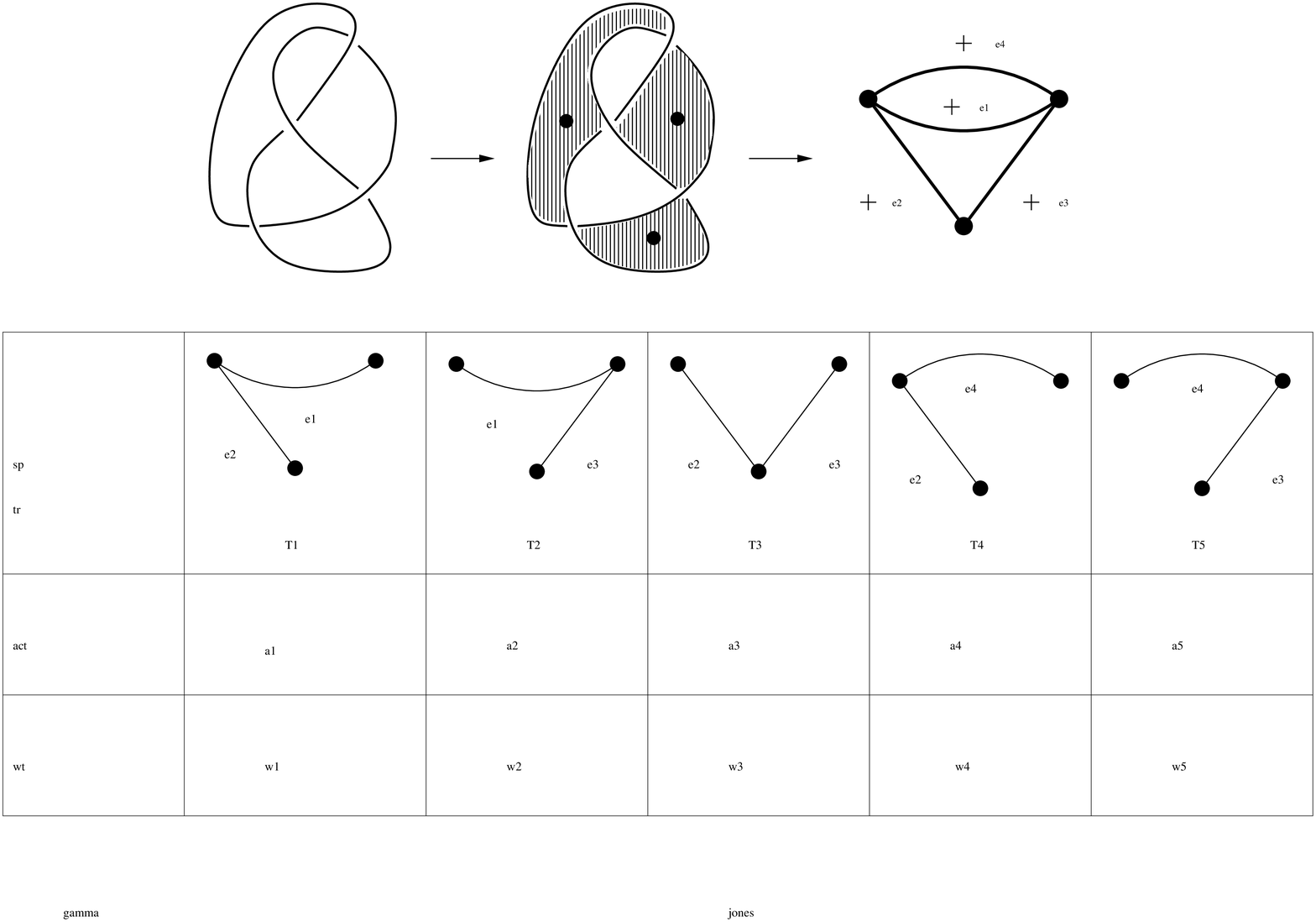} \\ \ \\
    \psfrag{A}{$A$}
    \psfrag{B}{$B$}
    \psfrag{U1}{\footnotesize{$U(T_1)$}} \psfrag{U2}{\footnotesize{$U(T_2)$}}
    \psfrag{U3}{\footnotesize{$U(T_3)$}} \psfrag{U4}{\footnotesize{$U(T_4)$}}
    \psfrag{U5}{\footnotesize{$U(T_5)$}}
    \includegraphics[height=4in]{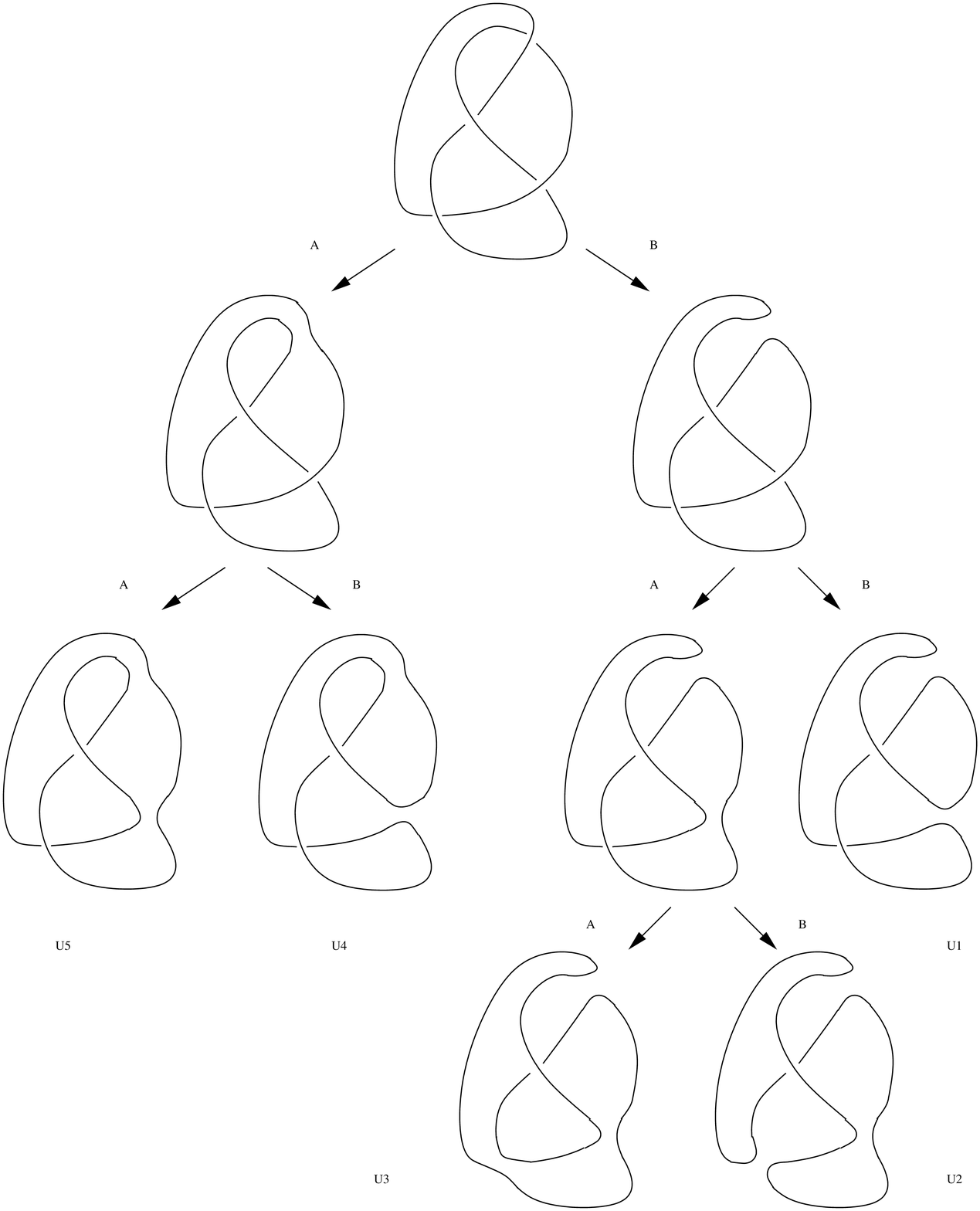}
  \end{center}
  \caption{Spanning trees and twisted unknots for figure-8 knot}
  \label{fig8}
\end{figure}

We can also consider each $U(T)$ as a partial smoothing of $D$
determined by $W(T)$.  In fact, there exists a skein resolution tree
for $D$ whose leaves are exactly all the partial resolutions $U(T)$,
for each spanning tree $T$ of $G$ (Theorem 2 \c{KHshort}).  Let
$\sigma(U) = \# A$-smoothings $- \#B$-smoothings, and let $w(U)$ be
the writhe.  If $U$ corresponds to $T$, then
$\mu(T)=A^{\sigma(U)}(-A)^{3w(U)}$ is exactly the monomial above Theorem
\ref{ThistThm}.  As Louis Kauffman pointed out, this is how humans
would compute $\kb{D}$: Instead of smoothing all the way to the final
Kauffman states, a human would stop upon reaching any twisted unknot
$U$, and use the formula $\mu(T)$.  We illustrate all of this for the
figure-eight knot diagram in Figure \ref{fig8}.

For any connected link diagram $D$, we choose the checkerboard
coloring such that its Tait graph $G$ has more positive edges than
negative edges, and in case of equality that the unbounded region is
unshaded.  In \c{KHshort}, we defined the spanning tree complex
$\C(D)=\{\C_v^u(D), \del\}$, whose generators correspond to spanning
trees $T$ of $G$.  The $u$ and $v$--grading are determined by $W(T)$
as follows:
$$ u(T)  =   \# L - \# \ell - \# \bar{L} + \# \bar{\ell} \quad {\rm and }\quad v(T) =  \# L + \# D $$
\begin{thm}[\c{KHshort}]\label{KHthm}
  For any connected link diagram $D$, there exist spanning tree
  complexes $\C(D)=\{\C_v^u(D), \del\}$ and $\UC(D)=\{\UC_v^u(D),
  \del\}$ with $\del$ of bi-degree $(-1,-1)$ that are deformation
  retracts of the reduced and unreduced Khovanov complexes,
  respectively.
\end{thm}

The differential in $\C(D)$ is defined indirectly.  As discussed in
detail below, for each $T$ the Khovanov complex $\widetilde{C}(U(T))$
is contractible, and we proceed by a sequence of collapses of each
$\widetilde{C}(U(T))$ to a single generator $Z(T)$.  The differential
on spanning trees is the one induced by all such collapses.  

Note that $u(T) = -w(U(T))$.  
Interestingly, $v(T)$ has appeared in several guises elsewhere:  
(1) Rasmussen's $\delta$--grading (Definition 4.4 \c{Rasmussen}) satisfies $\delta = 2v + k$, 
where $k$ is a constant that depends only on $D$.  
(2) A connected link diagram determines a {\em ribbon graph}, 
which is a graph embedded in a surface such that its complement is a union of $2$--cells.
The genus $g$ of the ribbon graph is the genus of the minimal such surface.
Each spanning tree of $G$ corresponds to a ribbon graph with one complementary $2$--cell, whose genus satisfies $g+v=k'$, 
where $k'$ is another constant that depends only on $D$ (Theorem 2.1 \c{dkh}).

\subsection{Fundamental cycle of a twisted unknot}\label{Sec_fundcyc}
We review the main ideas underlying Theorem \ref{KHthm}, which will be used in the next section.

Let $D$ be a connected link diagram with a basepoint $P$ away from the crossings of $D$.
In the version of Khovanov homology in \c{Viro, Viro2}, generators of the
reduced Khovanov complex $\widetilde{C}(D)$ are given by {\em enhanced
  Kauffman states} of $D$.  A Kauffman state $s$ is a choice of
smoothings of all crossings of $D$, and enhancements are $\pm$ signs
on every loop of $s$.  The reduced Khovanov complex consists only of
enhanced states for which every loop that contains $P$ has a positive
enhancement.  Enhanced states are incident in $\widetilde{C}(D)$ if
and only if exactly one $A$ marker can be changed to a $B$ marker,
such that loops unaffected by the marker change keep their
enhancements, and the changed loops are enhanced to increase the
enhancement signature by one.

\begin{figure}
\begin{center}
\psfrag{-}{$-$}
\psfrag{+}{$+$} \psfrag{=}{$=$}
\includegraphics[height=0.9in]{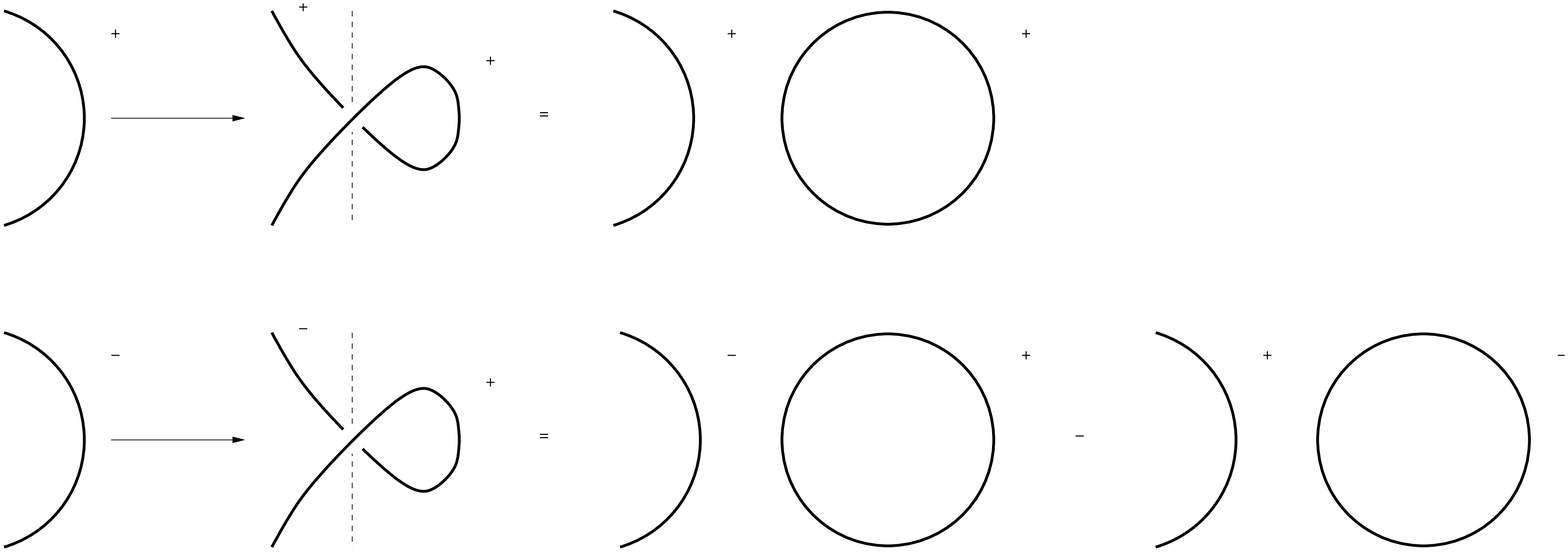}
\hspace*{1.3cm}
\includegraphics[height=0.9in]{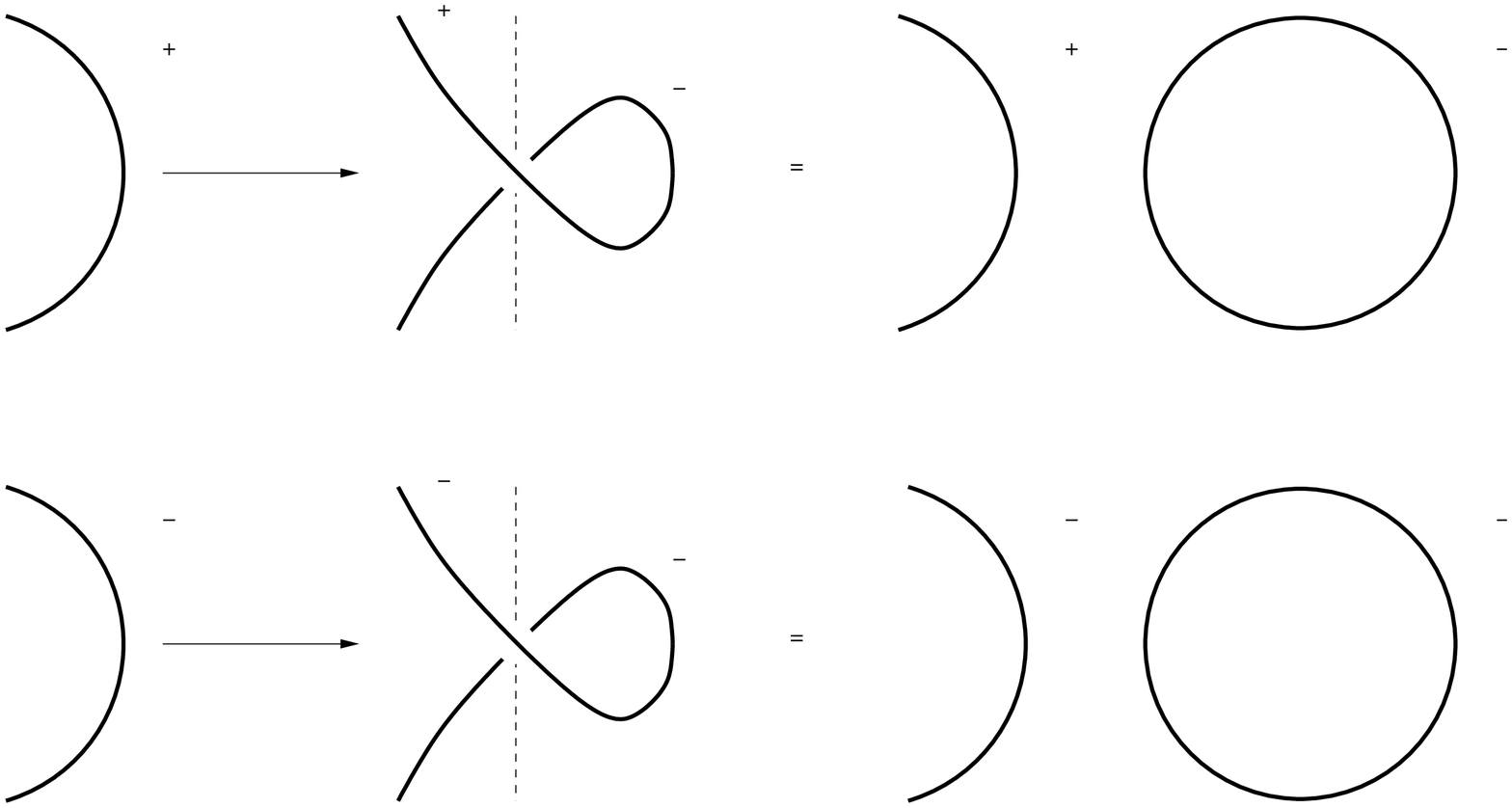}
\end{center}
\caption{How to obtain the fundamental cycle of a twisted unknot}
\label{twistfig}
\end{figure}

For any twisted unknot $U$, $\widetilde{C}(U)$ is contractible, with
the same homology as that of the positively enhanced round unknot
$\bigcirc^+$.  Starting from the round unknot, by a sequence of positive
and negative twists, we can obtain any $U$.  For every such twist,
Figure \ref{twistfig} indicates how to obtain a linear combination of
maximally disconnected enhanced states.  We define the {\em
  fundamental cycle} $Z_U \in \widetilde{C}(U)$ to be the linear
combination of enhanced states of $U$ given by iterating the local
changes in Figure \ref{twistfig}. Let $f_U: {\widetilde C}(\bigcirc)
\to {\widetilde C}(U)$ be defined by $f_U(\bigcirc^{+})=Z_U$.  

On the other hand, there exists a sequence of elementary collapses
$r_U: {\widetilde C}(U)\to {\widetilde C}(\bigcirc)$, such that
$r_U\circ f_U=id$ and $f_U\circ r_U\simeq id$.  Essentially, this
follows from invariance of Khovanov homology under the first
Reidemeister move \c{Khovanov}.  

We can summarize this discussion as follows: Each spanning-tree
generator $T\in \C(D)$ corresponds to a contractible Khovanov
subcomplex $\widetilde{C}(U(T))$, for which the fixed point of the
retraction is the fundamental cycle $Z_{U(T)}$.
The basepoint $P$ determines a {\em unique} fixed point for this
retraction: $(U,P)$ is given by a sequence of first Reidemeister moves
from $(\bigcirc,P)$, which determines $Z_U$ uniquely.

Let $\iota : {\widetilde C}(U(T)) \to {\widetilde C}(D)$ be the
inclusion of enhanced states given by appropriately shifting the
gradings.  The image $Z(T)=\iota\left(Z_{U(T)}\right)\in {\widetilde
  C}(D)$ is called the {\em fundamental cycle of} $T$.  
Note that $Z(T)$ is not generally a cycle in ${\widetilde C}(D)$, even 
though $Z_{U(T)}$ is a cycle in $\widetilde C(U(T))$. 
In the proof of Theorem \ref{KHthm},
the map from $\C(D)\to {\widetilde C}(D)$ given by $T\to Z(T)$ induces 
an isomorphism on homology.

Up to linear combinations of enhancements, $Z(T)$ is just a single
Kauffman state: the maximally disconnected state of $U(T)$, obtained
by replacing every positive or negative twist in $U(T)$ by an $A$ or
$B$ marker, respectively.  So from Table \ref{Table1} we obtain the
markers for $Z(T)$ from the activity word $W(T)$:
\begin{eqnarray}\label{livedeadAB}
\begin{tabular}{cccccccc}
$L$ & $D$ & $\ell$ & $d$ & $\bar{L}$ & $\bar{D}$ & $\bar{\ell}$ & $\bar{d}$ \\
\hline
$B$ & $A$ & $A$ & $B$ & $A$ & $B$ & $B$ & $A$ \\
\end{tabular}
\end{eqnarray}
It also follows that distinct enhanced states $s,s'\in
\widetilde{C}(U(T))$ differ only at markers that are live in $W(T)$.
If $i\neq j$, the enhanced states $s_i\in\widetilde{C}(U(T_i))$ and
$s_j\in\widetilde{C}(U(T_j))$ differ in at least one marker that is
dead in both $W(T_i)$ and $W(T_j)$.

Finally, it is straightforward to extend these ideas to unreduced
Khovanov homology.  Using the gradings in \c{Viro}, ${\widetilde
  C}(\bigcirc)\cong \Z^{(0,-1)}$ and
$C(\bigcirc)\cong\Z^{0,1}\oplus\Z^{0,-1}$.  Hence, for every $T$,
there are two fundamental cycles for $U(T)$, and two corresponding
generators in $\UC(D)$: $T^+$ in grading $(u(T),v(T))$, and $T^-$ in
grading $(u(T)+2,v(T)+1)$.  With the activity word $W(T)$ and the
basepoint $P$, we can associate a unique generator in $C(D)$ to each
of $T^+$ and $T^-$ by using the same rules in Figure \ref{twistfig} to
obtain $Z^{\pm}_{U(T)}\in C(D)$, starting with $\bigcirc^+$ for $T^+$
and $\bigcirc^-$ for $T^-$.

\subsection{Activity words and the differential on the spanning tree complex} \label{Sec_WD}

The proof of Theorem \ref{KHthm} does not provide an intrinsic
description of the differential on the spanning tree complex $\C(D)$
without reference to enhanced states.  The main result of this section is that the
simplest kind of incidence in $\C(D)$ is determined by activity words.

For a complex $(C,\del)$ over $\Z$ with graded basis $\{e_i\}$, let
$\kb{\cdot,\cdot}$ denote the inner product defined by
$\kb{e_i,e_j}=\delta_{ij}$.
We say $x$ is \emph{incident} to $y$ in $(C,\del)$ if
$\kb{\del x, y}\neq 0$ and their \emph{incidence number} is $\kb{\del x, y}$.

Let $T_1, T_2$ be spanning trees with fundamental cycles $Z_1,Z_2\in
{\widetilde C}(D)$.  We define $T_1$ and $T_2$ to be {\em directly
  incident} if $\kb{\del Z_1,Z_2}\neq 0$ in ${\widetilde C}(D)$.  In
this case, $\kb{\del Z_1,Z_2}=(-1)^{\beta}$, where $\beta$ is the
number of $B$--markers after the $A$--marker that is changed.  By
Lemma \ref{collapse4} below, if $T_1$ and $T_2$ are directly incident,
then they are incident in $\C(D)$ and $\kb{\del T_1,T_2}=\kb{\del
  Z_1,Z_2}=\pm 1$.
However, $T_1$ and $T_2$ may be incident in $\C(D)$ even though $\kb{\del Z_1,Z_2}=0$, which is discussed in Section \ref{Sec_higher}.  

\begin{thm}\label{direct}
Spanning trees $T_1$ and $T_2$ are directly incident if and only if
 the activity words $W(T_1)$ and $W(T_2)$ differ by changing exactly two
 (not necessarily adjacent) letters in one of the following four ways:
\begin{eqnarray*}
 L \; \bar{d} & \to &  d\;  \bar{D}  \\
 \bar{d}\;  D  & \to &  \bar{L}\;  d  \\
 \bar{\ell}\;  D  & \to &  \bar{D}\;  d  \\
 D\;  \bar{d}  & \to &  \ell\;  \bar{D}  \\
\end{eqnarray*}
In particular, $T_2$ is obtained from $T_1$ by replacing one positive
edge $e\in T_1$ with one negative edge $f$, such that $f \in
cut(T_1,e)$, and no other edges change activity.
\end{thm}
\pf
First, we show that if $W(T_1)$ (on the left) changes in one of the
four ways to $W(T_2)$, then $T_1$ and $T_2$ are directly incident.
Let $Z_1$ and $Z_2$ be fundamental cycles of $T_1$ and $T_2$.  In all
four cases, by (\ref{livedeadAB}) exactly one $A$ marker of $Z_1$ is
changed to a $B$ marker to get $Z_2$, and $(u(T_2), v(T_2)) = (u(T_1)-1,
v(T_1)-1)$.  Changing indices according to equations (2) in \c{KHshort}, it
follows by results in \c{Viro} that at least one summand of each of
$Z_1$ and $Z_2$ are incident in $\widetilde{C}(D)$.

We claim that $\kb{\del Z_1,Z_2}\neq 0$.  If these are single enhanced
states, then we are done.  For linear combinations of enhanced states,
we must show that incidences among summands do not cancel.  
A fundamental cycle $Z(T)$ can have more than one summand only if
$U(T)$ is smoothed at a crossing $c$, resulting in a linear
combination of enhanced states, as shown in Figure \ref{twistfig}.
Since $c$ is a crossing of $U(T)$, $c$ is live in $W(T)$.
In all four cases, the marker that changes from $A$ to $B$ is dead in both $W(T_1)$ and $W(T_2)$, so the marker at $c$ cannot change.
All summands of $Z(T)$ have the same markers, so the sign of every summand is determined by its enhancements.
Since the sign of the Khovanov differential depends only on the markers, 
cancellations cannot occur among terms in $\kb{\del Z_1,Z_2}$.
Since at least some summands of $Z_1$ and $Z_2$ are
incident and do not cancel, $T_1$ and $T_2$ are directly incident.

Conversely, suppose $T_1$ and $T_2$ are directly incident.  We claim
there is exactly one pair of edges $e_i,\, e_j$ such that $T_2 = (T_1
\setminus e_i)\cup e_j$, and only $e_i$ and $e_j$ change activities.

If a marker does not change, then by (\ref{livedeadAB}), since edge signs
do not change, the activity of the corresponding edge can change as
follows:
\begin{equation}\label{activityswitch}
L\leftrightarrow d,\qquad  D\leftrightarrow\ell,\qquad \bar{L}\leftrightarrow\bar{d},\qquad \bar{D}\leftrightarrow\bar{\ell}
\end{equation}
Therefore, without a marker change, the activity of an edge changes if
and only if the edge is removed from the tree or inserted into the
tree.

  From any spanning tree $T$, we can obtain any other spanning tree 
by switching pairs of edges $e_i\in T,\, e_j\notin T$, such that $e_j\in cut(T,e_i)$.
Consider switching one such pair of edges for which neither marker changes.

Suppose the markers of $e_i$ and $e_j$ are fixed, and suppose for
spanning trees $T,T'$, we have $T' = (T \setminus e_i)\cup e_j$.  In
every case in (\ref {activityswitch}), $e_i$ and $e_j$ are both live
in either $T$ or $T'$.  However, $e_j\in cut(T,e_i)$ and $e_i\in
cut(T',e_j)$, so only one of $e_i$ or $e_j$ can be live (the
lower-ordered edge).  This contradiction implies that if neither
marker changes, then the activities cannot change, and in particular,
this pair of edges cannot be switched.

Since $T_1$ and $T_2$ are directly incident, exactly one marker
changes.  By the argument above, there is exactly one pair of edges
$e_i,\, e_j$ such that $T_2 = (T_1 \setminus e_i)\cup e_j$, and only
the activities of $e_i$ and $e_j$ change.  Moreover, only the
lower-ordered edge can be live in either $T_1$ or $T_2$.  Since
$v(T_2)= v(T_1)-1$, $e_i$ must be positive, and $e_j$ negative.  Since
$u(T_2)= u(T_1)-1$, if both edges are dead on the right (i.e., with
respect to $T_2$), one edge on the left must be $L$ or $\bar{\ell}$;
if both edges are dead on the left, one edge on the right must be
$\bar{L}$ or $\ell$.  These four cases are the ones given in the
theorem, and all can occur.  \done

\begin{lemma}\label{collapse4}
Let $T_1, T_2$ be spanning trees with fundamental cycles $Z_1,Z_2\in {\widetilde C}(D)$.
If $\kb{\del Z_1,Z_2}\neq 0$ then in $\C(D)$, $\kb{\del T_1,T_2}=\kb{\del Z_1,Z_2}$.
\end{lemma}
\pf 
If $x$ is incident to $y$ in $\widetilde{C}(D)$, we denote this by $x\to y$ below.
Let $U_i = U(T_i)$.
We claim that the differential $Z_1 \to Z_2$ remains
after all elementary collapses of twisted unknots, as in Lemma 4 of
\c{KHshort}.  It suffices to show that the incidences shown in the diagram
below are impossible for any enhanced states $s',\, s''$ that are
distinct from $Z_1,\, Z_2$.  This is the only way for the differential
$Z_1 \to Z_2$ to be removed by elementary collapse.
$$
\xymatrix{
{s'}\ar[dr] \ar[r] & {s''} \\
{Z_1}\ar[ur] \ar[r] & {Z_2} }
$$

{\bf Case 1:}\quad $s'\in \widetilde{C}(U_1)\subset \widetilde{C}(D)$.
If $i\neq j$, any incidence between enhanced states in
$\widetilde{C}(U_i)$ and $\widetilde{C}(U_j)$ must occur at a marker
that is dead in both $W(T_i)$ and $W(T_j)$.  Thus, both $s'$ and $Z_1$
differ from $Z_2$ on a dead marker, hence they have the same live
markers.  Since both are in $\widetilde{C}(U_1)$, they have the same
dead markers too.  Therefore, $s'$ and $Z_1$ just differ by the
following enhancements:
$$
\xymatrix{
{\bigcirc^- \bigcirc^+}\ar[dr] \ar[r] & {s''} \\
{\bigcirc^+ \bigcirc^-}\ar[ur] \ar[r] & {\bigcirc^+} }
$$
Now, for both $s'$ and $Z_1$ to be incident to $s''$, the same marker
must change.  This implies that $s''=Z_2$ since both have the same
markers and the same enhancements, which is a contradiction.

{\bf Case 2:}\quad $s'\not\in \widetilde{C}(U_1)\subset \widetilde{C}(D)$.
Suppose $Z_1 \to  Z_2$ at marker 1, and $s' \to Z_2$ at marker 2, which must be
distinct markers.  Therefore at markers 1 and 2, we have
$$
\xymatrix{
{BA}\ar[dr] \ar[r] & {s''} \\
{AB}\ar[ur] \ar[r] & {BB} }
$$
Because $Z_1$ and $s'$ are both incident to $s''$, this implies that $s''$ must have the same markers as $Z_2$.
Therefore, for $Z_1$ to be incident to both $s''$ and $Z_2$, the enhancements
must be as follows:
$$
\xymatrix{
{s'}\ar[dr] \ar[r] & {\bigcirc^- \bigcirc^+} \\
{\bigcirc^-}\ar[ur] \ar[r] & {\bigcirc^+\bigcirc^-} }
$$
Now, for both $s''$ and $Z_2$ to be incident to $s'$, the same marker
must change.  So marker 1 = marker 2, which is a contradiction.
\done

\subsection{Spanning tree filtration and spectral sequence} \label{Sec_SS} 

The activity word $W(T)$ determines a partial smoothing $U(T)$.
Live edges, denoted below by $*$, correspond to crossings of the twisted unknot, which are not smoothed.

Let $D$ be any connected link diagram with $n$ ordered crossings.
For any spanning trees $T, T'$ of $G$, let $(x_1,\ldots,x_n)$ and $(y_1,\ldots,y_n)$ be the corresponding partial smoothings of $D$.
We define a relation $T > T'$, or equivalently,
$(x_1,\ldots,x_n) > (y_1,\ldots,y_n)$
if for each $i$, $y_i=A$ implies $x_i = A$ or $*$, and there exists $i$ such that $x_i=A$ and $y_i=B$.
The transitive closure of this relation gives a partial order, also denoted by $>$ (Proposition 1 \c{KHshort}).
We define $\mathcal P(D)$ to be the poset of spanning trees of $G$ with this partial order.
Note that $\mathcal P(D)$ always has a unique maximal tree and unique minimal
tree, whose partial smoothings contain the all-$A$ and all-$B$
Kauffman states, respectively.

For example, for the figure-8 knot from Figure \ref{fig8},
\begin{center}
\begin{tabular}{c|c|c|c|c}
$T_1$ & $T_2$ & $T_3$ & $T_4$ & $T_5$ \\
\hline
$**BB$ & $*BAB$ & $*AAB$ & $**BA$ & $**AA$ \\
\end{tabular}
\end{center}
We get two sequences:
$ T_5 > T_3 > T_2 > T_1$ and $T_5 >T_4>T_1$.
The maximal and minimal trees correspond respectively to the left-most and right-most unknots in Figure \ref{fig8}.

The poset $\mathcal P(D)$ provides a partially ordered filtration of
${\widetilde C}(D)$ indexed by $\mathcal P(D)$: 
Let $\U_i = \iota({\widetilde C}(U(T_i))) \subset {\widetilde C}(D)$.
Let $\psi: \mathcal P(D)\to {\widetilde
  C}(D)$ defined by $\psi(T) = +_{T \geq T_i}\U_i$.  From a partially ordered
filtration, we can get a decreasing linearly ordered filtration
$F^p{\widetilde C}(D)$ by taking all trees of order at least $p$ from
all maximal descending ordered sequences in $\mathcal P(D)$.
For example, the figure-8 knot from Figure \ref{fig8} has the following filtration $F^p {\widetilde C}(D)$,
$$ F^1 = \psi(T_5) = {\widetilde C}(D), \ \ F^2 = \psi(T_3) + \psi(T_4), \ \ F^3 = \psi(T_2), \ \ F^4 = \psi(T_1)$$

\begin{thm}[\c{KHshort}]\label{SS}
For any knot diagram $D$, there is a spectral sequence $E_r^{*,*}$
that converges to the reduced Khovanov homology
$\widetilde{H}^{*,*}(D;\Z)$, such that as groups $E^{*,*}_1 \cong
\C^*_*(D)$, and the spectral sequence collapses for $r\leq c(D)$, where
$c(D)$ is the number of crossings.
\end{thm}

The associated graded module consists of submodules of ${\widetilde C}(D)$ in bijection with spanning trees:
\begin{equation} \label{E0}
E^{p,*}_0 = F^p{\widetilde C}(D)/F^{p+1}{\widetilde C}(D) = \oplus_i\ \U_i 
\end{equation}

\begin{cor}\label{cor_E2W}
  For any knot diagram $D$, the $E_2$--term of the spectral sequence
  in Theorem \ref{SS} is determined by the set of activity words for all spanning trees in $G(D)$.
\end{cor}
\pf Let $D$ be any knot diagram.  It follows from the filtration that
for any $p$, if $\U_1, \U_2\subset E^{p,*}_0$, then $T_1$ and $T_2$
are not comparable in $P$.  Hence, $d_0: E^{p,q}_0 \to E^{p,q+1}_0$
satisfies $d_0(\U_k)\subset\U_k$ for every $k$. 
This implies that (\ref{E0}) is a direct sum of complexes. 

Each complex $\U_k$ has homology generator corresponding to a spanning tree,
so $E_1$ is isomorphic as a group to the spanning tree complex:
$$E^{*,*}_1 =H^*(F^p/F^{p+1},d_0)= \oplus_k\ H^*(\U_k) \cong \C(D)$$
Let $d_1:E_1^{p,q}\to E_1^{p+1,q}$.
If $T_1 > T_2$ are directly incident and one filtration level apart
 then $Z_1\in F^p {\widetilde C}(D),\ Z_2\in F^{p+1} {\widetilde C}(D) \in E_1^{*,*}$, and $\kb{\del(Z_1),Z_2} \neq 0$. 
Thus, $\kb{d_1(Z_1),Z_2} \neq 0$.
Conversely, if $\kb{d_1(Z_1),Z_2} \neq 0$ then $T_1$ and $T_2$ are 
one filtration level apart and hence directly incident.
The partial order
and filtration are determined by activity words, and by Theorem
\ref{direct} direct incidences are determined by activity words. 
Therefore, the $E_2$--term of the spectral sequence is determined 
by activity words. 
\done

\subsection{Higher order incidences} \label{Sec_higher}
To construct $\C(D)$ as well as $\UC(D)$, we proceed by a sequence of
elementary collapses of each $\U_i$ to its fundamental cycle $Z(T_i)$,
starting from the minimal tree and ascending in the partial order
whenever trees are comparable.  Any elementary collapse in $\U_i$ does
not change incidence numbers in $\U_j$ for any $j\neq i$ (Lemma 5
\c{KHshort}), so we can sequentially collapse each $\U_i$.

Differentials are induced from each collapse, so $T_1$ and $T_2$ may
be incident in $\C(D)$ without being directly incident; i.e.,
$\kb{\del Z_1,Z_2}=0$ but $\kb{d_r Z_1,Z_2}\neq 0$ for some $r>1$.  In
general, $\kb{\del T_1, T_2}$ is the sum of induced incidence numbers
given by all {\em ladders} from $Z_1$ to $Z_2$.  Before giving
definitions, here is an example of a $2$--incidence ($i\neq 1,2$ and $x_i,\, y_i\neq Z_i$), which becomes an incidence after collapsing $U_i$:
$$
\xymatrix{
{Z_1} \ar[dr]         &        &    U_1 \\
{x_i} \ar[dr] \ar[r]  & {y_i}  &    U_i \\
                      & {Z_2}  &    U_2
} $$
\begin{defn}\label{kdef}
Let $T_1 > T_2$ be spanning trees with fundamental cycles $Z_1,Z_2\in {\widetilde C}(D)$.
$T_1$ and $T_2$ are $1$--incident if they are directly incident.
For $k>1$, $T_1$ and $T_2$ are {\em $k$--incident} if
there exist $x_i, y_i \in \U_{j_i}-Z_{j_i}$ for $1\leq i\leq k-1$,
such that if $x_0=Z_1$ and $y_k=Z_2$ then
$$ \kb{\del x_i,y_{i+1}}\neq 0 \ {\rm for }\ 0\leq i\leq k-1, \quad \kb{\del x_i,y_i}\neq 0 \ {\rm for }\ 1\leq i\leq k-1 $$
Such a sequence of ordered pairs of enhanced states will be
called a {\em ladder of enhanced states} from $Z_1$ to $Z_2$.
\end{defn}

If $T_1$ and $T_2$ are $k$--incident, then collapsing along $\{\U_{j_i}\ |\ 1\leq i\leq k-1,\ j\neq 1,2\}$, 
as in Definition \ref{kdef}, the incidence number between $Z_1$ and $y_k=Z_2$ induced from this ladder is
$\displaystyle (-1)^{k-1}\kb{\del Z_1,y_1}\prod\nolimits_{i=1}^{k-1}\kb{\del x_i,y_i}\kb{\del x_i,y_{i+1}}$. 
Moreover, since each enhanced state belongs to a unique $\U_i$, for every ladder from $Z_1$ to $Z_2$, such a collapse implies the following:
\begin{enumerate}
\item $T_1 > T_{j_1} \geq \ldots \geq T_{j_{k-1}} > T_2$ and these relations are transitive.
\item Exactly $k$ $A$-markers of $Z_1$ are
changed to $B$-markers of $Z_2$.
\item Exactly $(k-1)$ $B$-markers of $Z_1$ are
changed to $A$-markers of $Z_2$.
\item For each $j_i,\ x_i,y_i$ are resolutions of $\U_{j_i}$, with a differential $x_i\to y_i$ by changing a marker that is live
  in $W(T_{j_i})$.
\item If at level $i$ and level $i+1$ ($i\neq 1,k-1$), the same spanning tree $T$
  occurs, then there exists a differential $x_i\to y_{i+1}$ by
  changing a marker that is live in $W(T)$.
\item If at level $i$ and level $i+1$, distinct spanning trees occur,
  then there exists a differential $x_i\to y_{i+1}$ by
  changing a marker that is dead in both $W(T_{j_i})$ and $W(T_{j_{i+1}})$.
\end{enumerate}

We will say that a sequence of (possibly repeated) activity words
$\{W(T_{j_i})\}$ along with the following extra data is an {\em
  admissible activity sequence} if
\begin{enumerate}
\item The first and last activity words correspond to spanning trees $T_1$ and $T_2$, whose bigradings permit a nonzero differential.
\item The sequence satisfies the partial order: $T_1 > T_{j_i} \geq T_{j_{i+1}} > T_2$ for all $i$.
\item A sequence of ordered pairs of markers indicates how to change the live markers in each $W(T_{j_i})$.
\end{enumerate}

Each ladder of enhanced states gives rise to a unique admissible activity sequence.  
The converse is an open question that is fundamental to understanding the differential 
on the spanning tree complex:
\begin{question}\label{ques_ladder}
  Which admissible activity sequences correspond to ladders of
  enhanced states?
\end{question}

Theorem \ref{direct} answers Question \ref{ques_ladder} in the simplest case.

As discussed in Section \ref{Sec_fundcyc}, given $W(T)$ and basepoint
$P$, we can compute $Z(T)\in \widetilde{C}(D)$ from $U(T)$.  In
particular, $W(T_1)$ and $W(T_2)$ completely determine $Z_1$ and
$Z_2$.  But starting with $Z_1$, and just specifying allowed marker
changes may not be sufficient to produce a ladder (or possibly a
linear combination of ladders) of enhanced states to $Z_2$.  The
difficulty inherent in Question \ref{ques_ladder} is whether the
enhancements on the states ``take care of themselves,'' or whether the
enhancements can obstruct the existence of a ladder, given a sequence
of allowed marker changes from $Z_1$ to $Z_2$.

For the unreduced spanning tree complex $\UC(D)$, an admissible
sequence must also include the signs of the spanning tree generators:
$\{\pm W(T_{j_i})\}$.  With the signed activity word $\pm W(T)$ and
the basepoint $P$, we can compute $Z^{\pm}$ for each generator
$T^{\pm}$.  Because any enhancement is allowed on the state with $P$,
it seems less likely, given an admissible signed activity sequence, that
enhancements can obstruct a ladder.

In Section \ref{KHmatroids}, we show how this is directly related to the mutation-invariance of Khovanov homology.

\section{Mutation and matroids} \label{Sec_MM}

A {\em marked} $2$-tangle is a $2$-tangle contained in a round ball
such that its four endpoints are equally spaced around the equator of
the boundary sphere, called a {\em Conway sphere}.  Let $L$ be a link
that contains a marked $2$-tangle $\tau$.  A {\em mutation} of $L$ is
the following operation: Remove the Conway sphere containing $\tau$,
rotate it by $\pi$ about one of its three coordinate axes, and glue it
back to form the link $L'$.

The same operation can be described for any planar diagram $D$ of $L$.
The projection of the Conway sphere is a {\em Conway circle} that
meets $D$ in four points, which are the endpoints of the marked
$2$-tangle diagram contained in the disc.  A mutation of $D$ is then
given by one of the three corresponding involutions of the disc.
Diagrams $D$ and $D'$ are called {\em mutants} if $D'$ can be obtained
from $D$ by a sequence of mutations.

\subsection{Tait graphs and mutation}
There are two choices for the checkerboard coloring of $D$,
and the resulting Tait graphs are the planar duals of each other.  The
projection of $D$ is the medial graph of $G$, and the signs on $G$
determine the crossings of $D$. This determines a one-one
correspondence between checkerboard-colored link diagrams and planar
embeddings of signed graphs.  
In order to study mutation using Tait graphs, we define two moves on graphs:

 \begin{figure}
 \begin{center}
 \includegraphics[width=4in]{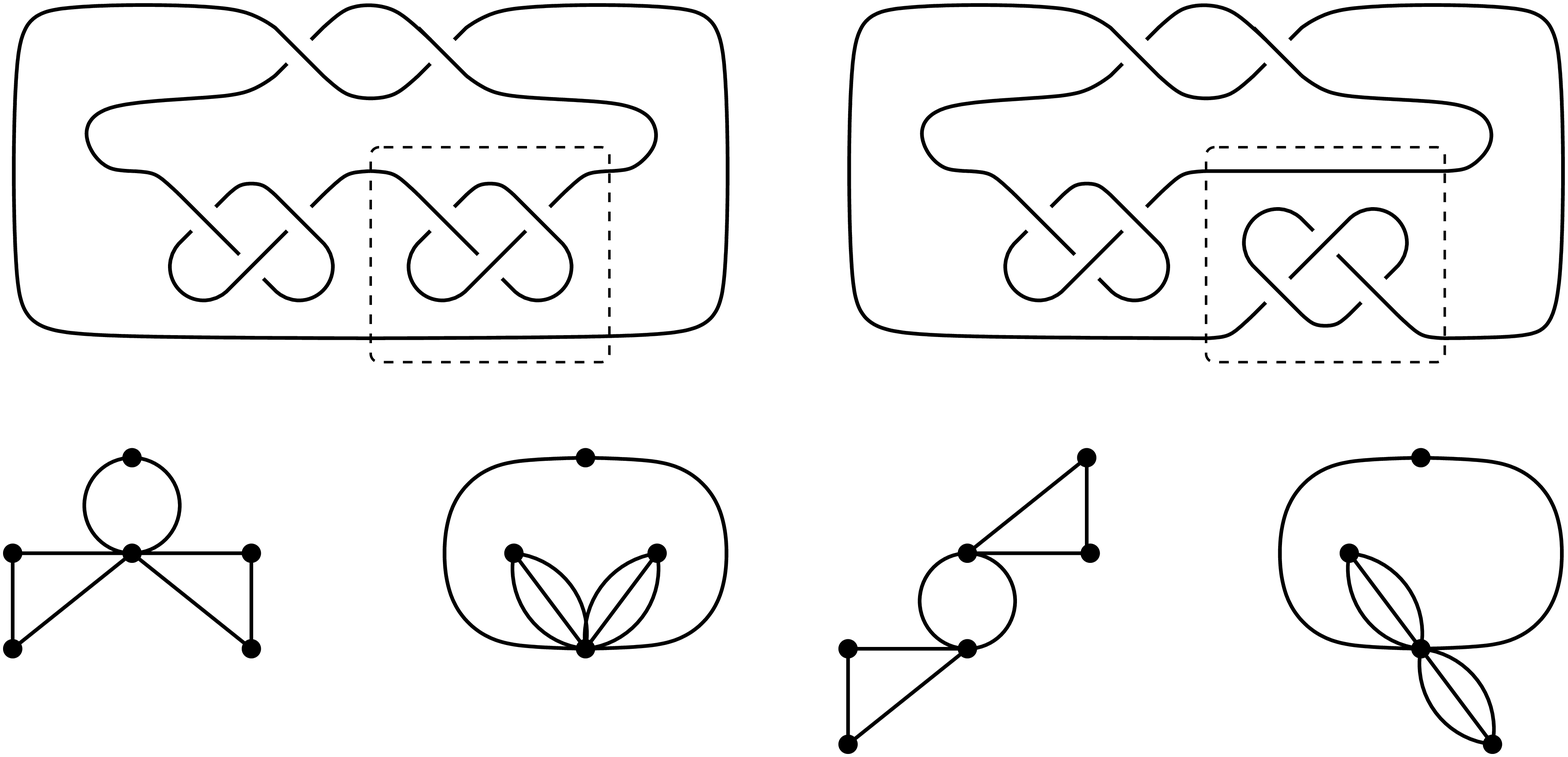}
\end{center}
 \caption{Connect sum for links and their Tait graphs}
 \label{connectsum}
 \end{figure}

{\bf $1$--flip}\quad  
Let $v_1$ and $v_2$ be vertices of disjoint graphs $G_1$ and
$G_2$.  A {\em vertex identification} is $G=G_1\sqcup G_2/v_1\sim v_2$.  If
$v$ is a cut-vertex of $G$, i.e. $G-v$ is disconnected, a {\em vertex
splitting} at $v$ of $G$ is the inverse operation of vertex
identification.  A $1$--flip of $G$ is a vertex splitting followed by a
vertex identification. 

{\bf $2$--flip}\quad 
For $i\in\{1,2$\}, let $u_i, v_i$ be vertices of disjoint graphs $G_i$ such 
that $G=G_1\sqcup G_2/(u_1,v_1)\sim (u_2,v_2)$.
A $2$--flip of  is the identification 
$G_1\sqcup G_2/(u_1,v_1)\sim (v_2,u_2)$.

We extend both these moves to signed graphs by requiring that 
the signs on the corresponding edges are preserved. 

For a link diagram $D$, a $1$--flip corresponds to breaking a connect sum 
and reconnecting at a different place. Since the connect sum operation 
is well-defined for knots, $1$--flips do not change the knot type. 
However, $1$--flips may change link type; see Figure \ref{connectsum}.
We will consider only component-preserving link mutation later.

\begin{figure}
 \begin{center}
 \includegraphics[width=4in]{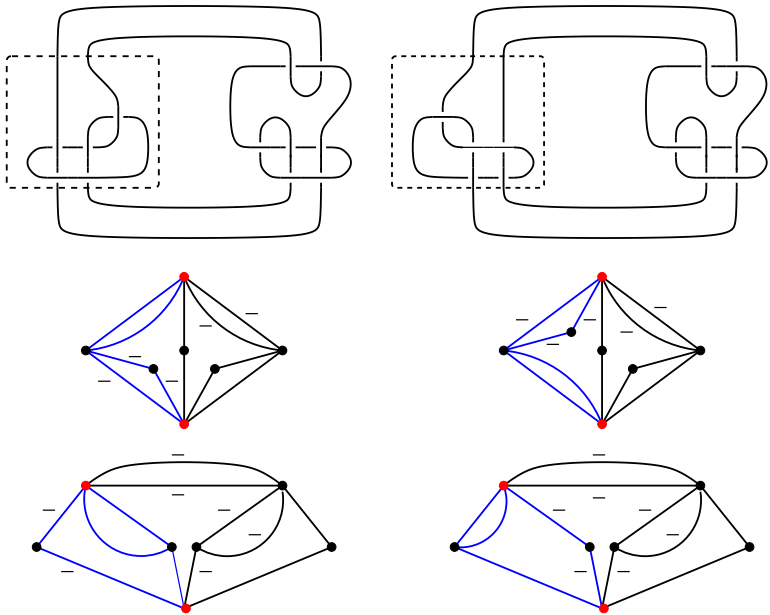}
\end{center}
 \caption{Kinoshita-Terasaka and Conway mutants and their Tait graphs}
 \label{kt}
 \end{figure}

$2$--flips correspond to mutation for link diagrams.  Figure \ref{kt}
shows the Kinoshita-Terasaka and Conway mutants along with their Tait
graphs (unsigned edges are positive). The graphs in the second row
come from the checkerboard coloring with the unbounded region shaded,
and the graphs in the third row from the other checkerboard coloring.

Some mutations change only the planar embedding of $G$ but not $G$
itself, so not all types of mutation can be realized as $2$--flips.
For example the graphs in the third row of Figure \ref{kt} are not
related by $2$--flips.  To address this, we define the following two
moves on planar embeddings of $G$ that preserve the graph itself.

{\bf planar $1$--flip}\quad 
A planar $1$--flip replaces a $1$-connected component of a planar embedding 
with its rotation by $\pi$ about an axis in the plane which intersects 
the cut vertex. 

{\bf planar $2$--flip}\quad 
A planar $2$--flip replaces a $2$-connected component of a planar embedding 
with its rotation by $\pi$ around the axis determined by the $2$-connecting
vertices. 

We similarly extend both these moves to embeddings of signed graphs by
requiring that the signs on the corresponding edges are preserved.

Any two planar embeddings of a signed graph are related by a sequence of
planar $1$--flips and planar $2$--flips (see \c{mt2001}).  
As before, these moves correspond to reconnecting connect sums and
mutations of link diagrams, respectively.  Although, $1$--flips can
also correspond to mutation in link diagrams whose Tait graphs have a
cut vertex; for example, see Figure \ref{connectsum} and \cite{Wehrli2}.

A graph $G$ is said to be $2${\em --isomorphic} to a graph $H$ if $G$
can be obtained from $H$ by any sequence of vertex identifications,
vertex splittings, or $2$--flips.  Hence, a connected graph $G$ is
$2$--isomorphic to a connected graph $H$ if $G$ can be obtained from
$H$ by any sequence of $1$--flips and $2$--flips.  In particular,
isomorphic graphs are $2$--isomorphic. We require 2--isomorphisms of 
  signed graphs to preserve signs on the edges. 
\ 
\begin{prop}\label{taitmutation}
Let $D$ and $D'$ be connected link diagrams with checkerboard
colorings chosen so that their unbounded regions are both shaded or
both unshaded.  Let $G$ and $G'$ be their respective Tait graphs.
Then $D$ and $D'$ are mutants if and only if $G$ and $G'$ are
$2$--isomorphic.
\end{prop}
\pf
For any Tait graph, any type of mutation corresponds to
either a $1$--flip (possibly a planar $1$--flip), a $2$--flip or a
planar $2$--flip, and all of these can be realized by mutation.
As mentioned above,
any two planar embeddings of a graph are related by a sequence of
planar $1$--flips and planar $2$--flips. 
Specifying the coloring of the unbounded region distinguishes a Tait graph from its planar dual.
\done

Thus, in order to study mutation via Tait graphs, we need to study
invariants of $2$--isomorphism classes of signed graphs.  As we discuss
below, these naturally come from matroids.

\subsection{Matroids}
We recall some ideas from the theory of matroids (see \c{Oxley}).  
A matroid $M$ is a finite set of elements, together with a family of subsets, called {\em independent sets}, such that
\begin{enumerate}
\item The empty set is independent,
\item Every subset of an independent set is independent,
\item For every subset $A$ of $M$, all maximal independent sets contained in $A$ have the same number of elements.
\end{enumerate}
A maximal independent set in $M$ is called a {\em basis} for $M$, and
any two bases of $M$ have the same number of elements, which is the
rank of M.

For example, let $E$ be the set of edges of a graph $G$, and let
$\mathcal{I}$ be the collection of subsets of edges that do not
contain a cycle.  Then $(E,\mathcal{I})$ is a matroid $M(G)$, called
the {\em graphic matroid} of $G$.  For a connected graph $G$, the
bases of $M(G)$ are the spanning trees of $G$.

For any connected link diagram $D$ with a checkerboard coloring and
 Tait graph $G$, let the {\em colored graphic matroid} $M(D)$ be 
the graphic matroid $M(G)$ with edges
colored by $\{\pm 1\}$ as in the Tait graph, according to the
crossings of $D$.  

Whitney \c{Whitney} determined precisely when two graphs have
isomorphic graphic matroids.  This fundamental result, which motivated matroid
theory, is called the $2$--isomorphism theorem (for background
see \c{Oxley}).  If we require that any isomorphism of colored graphic
matroids be color-preserving, then the $2$--isomorphism theorem
extends to signed graphs (see e.g., \c{VertiganWhittle}):

\begin{thm}
\label{WhitneyThm}
For signed graphs $G$ and $H$ with no isolated vertices, 
their colored graphic matroids are isomorphic if and 
only if $G$ and $H$ are $2$--isomorphic.
\end{thm}
 
Theorem \ref{WhitneyThm} and Proposition \ref{taitmutation} imply the following:
\begin{cor}\label{mut_cor}
Let $D$ and $D'$ be connected link diagrams with checkerboard
colorings chosen so that their unbounded regions are both shaded or
both unshaded.  Let $M(D)$ and $M(D')$ be their respective colored graphic matroids.
Then $D$ and $D'$ are mutants if and only if $M(D)\cong M(D')$.
\end{cor}

Consequently, any knot invariant $\varphi$ is invariant under mutation
if and only if for any knot diagram $K$, $\varphi(K)$ is an invariant
of the colored graphic matroid $M(K)$.  For any matroid $M$,
activities can be defined with respect to its basis, just as we did
for the graphic matroid $M(G)$ using spanning trees.  We will use that
activity words are determined by $M(K)$, essentially due to Crapo
\c{Crapo}.

For example, by Theorem \ref{ThistThm}, the Jones polynomial $V_K(t)$
has an expansion using activity words.  Therefore, the Jones
polynomial is an invariant of $M(K)$, and hence invariant under
mutation.  Below, we extend this idea to Khovanov homology.

\subsection{Khovanov homology and matroids}\label{KHmatroids}
For a given connected link diagram $D$ with basepoint $P$, we choose the checkerboard coloring such
that its Tait graph $G$ has more positive edges than negative edges,
and in case of equality that the unbounded region is unshaded.  Let
$M(D)$ be the colored graphic matroid of $D$ with this coloring.
The generators of $\C(D)$, which are the spanning trees of $G$, are
bases of $M(D)$.  Since both the $u$ and $v$--gradings are determined
by the activities and signs, the bigrading on $\C(D)$ is determined
by $M(D)$.

Generally, Conway mutation, as in Figure \ref{connectsum}, may not
preserve components.  Indeed, two mutant links were shown to have
different Khovanov homology in \c{Wehrli2}, using this connect sum
ambiguity for links.  From our point of view, such a mutation moves
the basepoint $P$ from one component to another, leading to a
different $Z_{U(T)}$ for every $T$, which sometimes changes the
homology.  To eliminate this ambiguity, we can either consider only
knot diagrams, or require that Conway mutation of links be
component-preserving.  For purposes of exposition, it is easier to
just discuss mutation of knot diagrams.

Whenever $K$ and $K'$ are mutant knot diagrams, by Corollary 
\ref{mut_cor}, $M(K)\cong M(K')$.  
Therefore, $\C(K) \cong \C(K')$ as bigraded abelian groups. We conjecture 
that the differential on $\C(K)$ is determined by $M(K)$ in the following way.

\begin{conj}\label{mut_conj}
Let $K$ and $K'$ be knot diagrams such that $M(K)\cong M(K')$. 
If $T_1, T_2\in\C(K)$ and $T_1', T_2'\in\C(K')$ are generators corresponding to spanning trees,  
$$ \kb{\del T_1, T_2} = \kb{\del T_1', T_2'} \quad {\rm whenever }\quad W(T_1)=W(T_1'),\ W(T_2)=W(T_2')$$
\end{conj}

If Conjecture \ref{mut_conj} holds, then $\C(K) \cong \C(K')$ as 
bigraded chain complexes for mutant knot diagrams $K$ and $K'$. 
This would imply that $\widetilde{H}(K)$ is invariant under mutation.

A quasi-isomorphism between chain complexes is a morphism that induces
an isomorphism on homology. Any two chain complexes of free abelian
groups with isomorphic homology are quasi-isomorphic.  
\footnote{This follows from the fact that every chain complex of free
  abelian groups decomposes as a direct sum of two-step complexes, for
  which the relation matrix can be diagonalized. We thank Ciprian
  Manolescu for this comment. } 
This implies the following equivalence:

{\em For a knot diagram $K$, the reduced Khovanov homology 
$\widetilde{H}(K)$ is invariant under mutation if and only if 
$\C(K)$ is determined up to quasi-isomorphism by $M(K)$. }

\begin{cor}\label{cor_E2M}
For any knot diagram, the $E_2$--term of the spectral sequence in Theorem \ref{SS} is invariant under mutation.
\end{cor}
\pf 
Corollary \ref{cor_E2W} implies that $E^{*,*}_2(K)$ is determined by $M(K)$.  If $K$ and $K'$ are mutant
knot diagrams, $M(K)\cong M(K')$, which implies $E^{*,*}_2(K)\cong
E^{*,*}_2(K')$.  \done

Theorem \ref{WhitneyThm} is at the heart of these results in terms of
spanning trees.  But in terms of enhanced Kauffman states, mutation
appears to be a rather violent operation on the Khovanov complex.  It
is interesting to relate these two points of view.

Generally, we are given a connected link diagram $D$ with a basepoint
$P$.  Conway mutation $\tau$ on $D$ induces a mutation on the Kauffman
states of $D$.  Conway mutation of an enhanced state $S$ of $D$ may
identify arcs with opposite enhancements.  To assign enhancements
unambiguously for $\tau(S)$, (1) any state disjoint from the Conway
circle must keep its enhancement, and (2) all enhancements must be
preserved in the part of $D$ that contains the basepoint $P$.  The
latter requirement induces enhancements on arcs in the other part of
$D$ that intersect the Conway circle.  We will refer to this
operation, which must preserve the link components, as Conway mutation
of $(D,P)$, denoted by $\tau(D,P)$.

Let $(D',P)=\tau(D,P)$.  By Theorem \ref{WhitneyThm}, spanning trees
$T' = \tau(T)$ if and only if $W(T)=W(T')$.  As discussed in Section
\ref{Sec_fundcyc}, with the activity word $W(T)$ and the basepoint
$P$, we can associate a unique generator in $C(D)$ to each of $T^+$
and $T^-,$ and for $\widetilde{C}(D)$ just use $T^+$.
  Thus, the activity word $\pm W(T_i)$ for $T_i^{\pm}$
determines a unique generator $Z_i^{\pm}$ in each of the respective
Khovanov complexes, $C(D)$ and $C(D')$, as well as in $\widetilde{C}(D)$ and $\widetilde{C}(D')$.  
Hence, the maps induced by
Conway mutation $\tau: \C(D)\to \C(D')$ and $\tau_{\mathcal U}: \UC(D)\to \UC(D')$ are isomorphisms of
bigraded abelian groups.

This provides an approach to prove mutation-invariance of Khovanov
homology.  In Section \ref{Sec_higher}, we defined an
admissible activity sequence, which depends only on $M(D)$.  Even
without an explicit answer to Question \ref{ques_ladder}, these
sequences may record the essential information about the differential:
\begin{conj}\label{conj_ladder}
  Let $(D',P)=\tau(D,P)$.  For every ladder of enhanced states in
  $\widetilde{C}(D)$, the corresponding admissible activity sequence describes a
  ladder of enhanced states in $\widetilde{C}(D')$.
\end{conj}
Conjecture \ref{conj_ladder} appears weaker than Conjecture
\ref{mut_conj}, but it too implies the mutation-invariance of reduced Khovanov
homology!

If Conjecture \ref{conj_ladder} holds then every differential in
$\C(D')$ may be computed from some collection of
admissible activity sequences.  If so, for every ladder in $\widetilde{C}(D)$,
there is a corresponding ladder in $\widetilde{C}(D')$ with the same induced
incidence number as in $\widetilde{C}(D)$.  This would imply that $\tau: \C(D)\to
\C(D')$ is a quasi-isomorphism.  In other words, 
$\C(D)$ is determined up to quasi-isomorphism by $M(D)$.

Because signs on the spanning trees (or their activity words) are not
contained in $M(D)$, the unreduced Khovanov complex $\UC(D)$ in
general is not determined by $M(D)$.  However, the following analogue
of Conjecture \ref{conj_ladder} similarly implies the mutation-invariance of unreduced Khovanov
homology:
\begin{conj}\label{conj_ladder_unred}
  Let $(D',P)=\tau(D,P)$.  For every ladder of enhanced states in
  $C(D)$, the corresponding admissible signed activity sequence describes a
  ladder of enhanced states in $C(D')$.
\end{conj}

\begin{rmk}\rm
  In \c{BNweb}, an attempt to prove mutation-invariance of Khovanov
  homology was outlined using ``re-embedding universality.''  However,
  as explained there, re-embedding universality implies invariance
  under cabled mutation, which Khovanov homology does not satisfy
  \c{DGST}.  We can explain the non-invariance of Khovanov homology
  under cabled mutation by the fact that cabled mutation corresponds
  to an $n$--flip for $n>2$.  Under this operation, Corollary
  \ref{mut_cor} 
  does not hold.  In fact, the $14$--crossing example in \c{DGST} has
  spanning trees whose activity words change after $2$--cabled
  mutation.  This lends some support to our approach.
\end{rmk}

\subsection*{Acknowledgments}
We thank Adam Sikora for helpful discussions, and the anonymous referee for thoughtful revisions. 

\bibliography{mkh}

\begin{thebibliography}{10}

\bibitem{BNweb}
D.~Bar-Natan.
\newblock {Mutation Invariance of {K}hovanov Homology},
  {http://katlas.math.toronto.edu/drorbn/}.

\bibitem{KHshort}
A.~Champanerkar and I.~Kofman.
\newblock Spanning trees and {K}hovanov homology.
\newblock arXiv:math.GT/0607510v3, to appear in {\em Proc. Amer. Math. Soc.}

\bibitem{dkh}
A.~Champanerkar, I.~Kofman, and N.~Stoltzfus.
\newblock Graphs on surfaces and {K}hovanov homology.
\newblock {\em Algebr. Geom. Topol.}, 7:1531--1540, 2007.

\bibitem{Crapo}
H.~Crapo.
\newblock The {T}utte polynomial.
\newblock {\em Aequationes Math.}, 3:211--229, 1969.

\bibitem{DGST}
N.~Dunfield, S.~Garoufalidis, A.~Shumakovitch, and M.~Thistlethwaite.
\newblock {Behavior of knot invariants under genus 2 mutation},
  arXiv:math.GT/0607258.

\bibitem{Rong_HG}
L.~Helme-Guizon and Y.~Rong.
\newblock A categorification for the chromatic polynomial.
\newblock {\em Algebr. Geom. Topol.}, 5:1365--1388, 2005.

\bibitem{Khovanov}
M.~Khovanov.
\newblock A categorification of the {J}ones polynomial.
\newblock {\em Duke Math. J.}, 101(3):359--426, 2000.

\bibitem{KhPatterns}
M.~Khovanov.
\newblock Patterns in knot cohomology. {I}.
\newblock {\em Experiment. Math.}, 12(3):365--374, 2003.

\bibitem{Kh_ICM2006}
M.~Khovanov.
\newblock Link homology and categorification.
\newblock In {\em International Congress of Mathematicians. Vol. II}, pages
  989--999. Eur. Math. Soc., Z\"urich, 2006.

\bibitem{LoeblMoffatt}
M.~Loebl and I.~Moffatt.
\newblock The chromatic polynomial of fatgraphs and its categorification.
\newblock {\em Adv. Math.}, 217(4):1558--1587, 2008.

\bibitem{MOS}
C.~Manolescu, P.~Ozsv{\'a}th, Z.~Szab{\'o}, and D.~Thurston.
\newblock On combinatorial link {F}loer homology.
\newblock {\em Geom. Topol.}, 11:2339--2412, 2007.

\bibitem{mt2001}
B.~Mohar and C.~Thomassen.
\newblock {\em Graphs on surfaces}.
\newblock Johns Hopkins Studies in the Mathematical Sciences. Johns Hopkins
  University Press, Baltimore, MD, 2001.

\bibitem{Oxley}
J.~Oxley.
\newblock {\em Matroid theory}.
\newblock Oxford University Press, 1992.

\bibitem{OZ1}
P.~Ozsv{\'a}th and Z.~Szab{\'o}.
\newblock Holomorphic disks and knot invariants.
\newblock {\em Adv. Math.}, 186(1):58--116, 2004.

\bibitem{Rasmussen}
J.~Rasmussen.
\newblock Knot polynomials and knot homologies.
\newblock In {\em Geometry and topology of manifolds}, volume~47 of {\em Fields
  Inst. Commun.}, pages 261--280. Amer. Math. Soc., Providence, RI, 2005.

\bibitem{thistlethwaite}
M.~Thistlethwaite.
\newblock A spanning tree expansion of the {J}ones polynomial.
\newblock {\em Topology}, 26:297--309, 1987.

\bibitem{VertiganWhittle}
D.~Vertigan and G.~Whittle.
\newblock A {$2$}-isomorphism theorem for hypergraphs.
\newblock {\em J. Combin. Theory Ser. B}, 71(2):215--230, 1997.

\bibitem{Viro2}
O.~Viro.
\newblock Khovanov homology, its definitions and ramifications.
\newblock {\em Fund. Math.}, 184:317--342, 2004.

\bibitem{Viro}
O.~Viro.
\newblock {Remarks on definition of {K}hovanov homology},
  arXiv:math.GT/0202199.

\bibitem{Wehrli2}
S.~Wehrli.
\newblock {Khovanov homology and {C}onway mutation}, arXiv:math.GT/0301312.

\bibitem{Wehrli3}
S.~Wehrli.
\newblock {Mutation invariance of Khovanov homology over $\Z/2\Z$}, talks at
  Knots in Washington, April 2007, and Kyoto University, May 2007.

\bibitem{Whitney}
H.~Whitney.
\newblock 2-{I}somorphic {G}raphs.
\newblock {\em Amer. J. Math.}, 55:245--254, 1933.

\end{thebibliography}
\bibliographystyle{plain}
\end{document}